\documentclass[a4paper, 12pt]{article}
\usepackage{amssymb, amsmath, amsfonts}
\usepackage[english,russian]{babel}
\usepackage{amsthm,amscd,amsfonts,latexsym,color,epsfig}

\begin{document}

\numberwithin{equation}{section}

\begin{center}
{\bf FRACTIONAL-DIFFRACTION-OPTICS
 \\ [3pt] CAUCHY PROBLEM:  \\ [3pt]
RESOLVENT-FUNCTION SOLUTION \\ [3pt] OF THE MATRIX INTEGRAL EQUATION}
\end{center}
\begin{center}
{\bf  Murat O. Mamchuev $^1$, Felix N. Chukhovskii $^1$}
\end{center}


\small

{\bf Abstract.} 
The fractional diffraction optics theory has been elaborated using the Green function technique. 
The optics-fractional equation describing the diffraction X-ray scattering by imperfect crystals has been derived as the fractional matrix integral Fredholm--Volterra equation of the second kind. 
In the paper, to solve the Cauchy problems, the Liouville--Neumann-type series formalism has been used to build up the matrix Resolvent-function solution. 
In the case when the imperfect crystal-lattice elastic displacement field is the linear function $f({\bf R}) = a x+b$, $a, b = const,$  the explicit solution of the diffraction-optics Cauchy problem has been obtained and analyzed for arbitrary fractional-order-parameter $\alpha$, $\alpha\in (0, 1].$
 \medskip

                  
                  Secondary 35F40, 35A08, 35C15, 35L40, 45F05, 35Q70, 35Q92;

 \smallskip


\normalsize

\section{Introduction}\label{sec:1}
\setcounter{equation}{0}

$\quad$ Usually, the diﬀraction optics (DO) has been formulated and based on the known differential Takagi--Taupin (TT) equations when the fractional-order-parameter $\alpha=1$ (see, e.g., 
 \cite{Takagi-ac-1969}, 
\cite{Takagi-1969}, \cite{Taupin-1964}, \cite{Bowen-1998}, \cite{Authier-2001}, \cite{Chukhovskii-1977}, \cite{Honkanen-2018},
\cite{Chukhovskii-2020} for details). 
At the same time, in the last decades, substantial progress has been achieved in mathematical physics using equations with fractional-order derivatives 
\cite{Oldham-1974}, \cite{Miller-Ross-1993}, \cite{Podlubny-1999}, \cite{Hilfer-2000}, \cite{n}, \cite{Kilbas-2006}, \cite{Uchaikin-2008}, \cite{Atanackovic-2014}.
Indeed, the different theory-physics models for solving the Cauchy problems using systems of fractional differential equations have been treated in
\cite{mamchuev-2010-du}, \cite{Mamchuev-2012},
\cite{Heibig-2012},
\cite{Kochubei-2013}, \cite{mam-mon-2013}, \cite{mamchuev-2015-mn}, 
\cite{mamchuev-2016-du}, 
\cite{mamchuev-2019-mnsvfu}, 
\cite{mamchuev-2021-aman}.

Following this logic, one can push one step further in the DO now founded on the TT-type equations with the fractional derivatives of the arbitrary order $\alpha\in (0, 1]$ along the crystal depth. 

In the paper, using Green’s function technique, one derives the matrix Fredholm–Volterra
integral equation of the second kind and builds up the matrix Resolvent-function solution of the DO
Cauchy problem. A goal of working is to develop the integral formalism of the DO theory earlier
proposed by authors \cite{MamChukh-2023} and to build up the matrix Resovent-function solution of the fractional Cauchy problem.

As an example, for an arbitrary fractional-order-parameter (FOP) $\alpha$, $\alpha\in (0, 1],$ when the crystal-lattice displacement field function $f({\bf R})$ is a linear function, namely: $f({\bf R}) = a x+b$, $a, b = const,$ one finds out an explicit analytical solution of the DO Cauchy problem.

Accordingly, the original system of fractional DO equations has the form (cf. \cite{MamChukh-2023})
$$\left (
\begin{array}{cc}
\partial_{0t}^{\alpha}-\frac{\partial}{\partial x}	& 0  	 \\
 0   	&      \partial_{0t}^{\alpha}+\frac{\partial}{\partial x}
\end{array}
\right )  
\left (
\begin{array}{c}
E_0(x,t)	 \\
E_h(x,t)
\end{array}
\right ) =$$
\begin{equation}\label{eq11}
=i
\left (
\begin{array}{cc}
\gamma						& \sigma \exp[if(x,t)]  	 \\
\sigma\exp[-if(x,t)]    	& \gamma	     
\end{array}
\right )  
\left (
\begin{array}{c}
E_0(x,t)	 \\
E_h(x,t)
\end{array}
\right ),
\end{equation}
with the DO Cauchy problem's condition
\begin{equation}\label{eq12}
\left (
\begin{array}{c}
E_0(x,0)	 \\
E_h(x,0)
\end{array}
\right ) =
\left (
\begin{array}{c}
\varphi_0(x)	 \\
\varphi_h(x)
\end{array}
\right ), 
\quad -\infty<x<\infty,
\end{equation}
where $\varphi_0(x)$ and $\varphi_h(x)$ are the given real-valued functions.

\section{Preliminaries}
\label{sec:2}
\setcounter{equation}{0}

Following up to \cite{n}, the Gerasimov--Caputo fractional derivative beginning at point $a$
is determined as (cf. \cite{n})
\begin{equation}\label{eqgc}
\partial_{at}^{\nu}g(t)={\rm sgn}^n(t-a)D_{at}^{\nu-n}\frac{d^n}{dt^n}g(t),
\quad n-1<\nu\leq n, \quad  n\in{\mathbb{N}},
\end{equation}
where
$D_{ay}^{\nu}$ is the Riemann--Liouville fractional integro-differential operator
of order $\nu$ is equal to  
$$D_{ay}^{\nu}g(y)=\frac{\mathop{\rm sgn}(y-a)}{\Gamma(-\nu)}
\int\limits_a^y \frac{g(s)ds}{|y-s|^{\nu+1}}, \quad  \nu< 0,$$
for $\nu \geq 0$ the operator $D_{ay}^{\nu}$ can be determined by 
recursive relation
\begin{equation}\label{eqrl}
D_{ay}^{\nu}g(y)=\mathop{\rm sgn}(y-a)\frac{d}{dy}D_{ay}^{\nu-1}g(y), \quad
\nu \geq 0,
\end{equation}
$\Gamma(z)$ is the Euler gamma-function.

Note that in the limit case of the FOP $\alpha=1$ the operator $\partial_{0t}^{\alpha}g(t)$ reduces to the standard  derivative 
$\frac{d}{dt}g(t).$

Formula for the composition of operators of fractional integration is valid as \cite{n}
\begin{equation}\label{eq21}
D_{0t}^{\nu}D_{0t}^{\delta}g(t)=D_{0t}^{\nu+\delta}g(t),
\quad \nu<0, \quad \delta<0.
\end{equation}
There is a formula for fractional integration by parts, namely \cite{n}:
\begin{equation}\label{eq22}
\int\limits_0^t g(t,\xi)D_{0\xi}^{\nu}h(\xi)d\xi=
\int\limits_0^t h(\xi)D_{t\xi}^{\nu}g(t,\xi)d\xi,
\quad \nu<0.
\end{equation}

Further, we denote the Fourier transform of the function $f(x)$ by $(f(x))_k$,  the Laplace transform of the function $g(t)$ by $(g(t))_p$,  and  respectively, the double Fourier--Laplace transform of the function $h(x,t)$ by
$(h(x,t))_{k,p}$. 

By using the following formula
\cite[p. 98]{Kilbas-2006}
$$[\partial_{0t}^{\alpha}H(x,t)]_{p}=p^{\alpha}[H(x,t)]_{p}-
p^{\alpha-1}H(x,0),$$ 
one can get
\begin{equation}\label{eq23}
[O_{\pm}^{\alpha}H(x,t)]_{k,p}=(p^{\alpha}\pm ik)[H(x,t)]_{k,p}-
p^{\alpha-1}[H(x,0)]_{k},
\end{equation}

The following series
\begin{equation*}
\phi(\beta, \rho;z)=\sum\limits_{k=0}^{\infty}\frac{z^k}{k!\Gamma(\beta k+\rho)},
\quad \beta>-1, \quad \rho\in {\mathbb C}
\end{equation*}
defines the Wright function
\cite{Wright-1934} which  depends on two parameters $\rho$ and $\mu.$


Accordingly, the following differentiation formula is valid  \cite{Wright-1934}:
\begin{equation}\label{eq24}
\frac{d}{dz}\phi(\beta, \rho; z)=\phi(\beta, \beta+\rho; z), \quad \beta>-1.
\end{equation}
Let $\beta\in (0,1),$ and $\mu, \nu \in {\mathbb R},$
and the inequality takes place
$$
\beta\in(0,1), \quad
0\leq |\arg \lambda|<\frac{1-\beta}{2}\pi,
$$
then the formula 
\begin{equation}\label{eq26}
D_{0y}^{\nu}y^{\mu-1}\phi(-\beta, \mu; -\lambda y^{-\beta})=y^{\mu-\nu-1}\phi(-\beta, \mu-\nu;  -\lambda y^{-\beta})
\end{equation}
is valid \cite{Pskhu-2019}.

The following formula  
take place \cite{Stankovic-1970}
\begin{equation}\label{eq27}
\left(y^{\delta-1}\phi(-\beta,\mu;-ty^{-\beta})\right)_p=p^{-\mu}e^{-p^{\beta} t}.
\end{equation}

Correspondingly, the integrals 
\begin{equation}\label{eq28}
\int\limits_{0}^{\infty}\frac{\cos kx}{k^2+\rho^2}dk=\frac{\pi}{2\rho}e^{-\rho x},
\end{equation}
\begin{equation}\label{eq29}
\frac{1}{2\pi i}\int\limits_{-i\infty}^{i\infty}
\frac{e^{-x\sqrt{p^2+\sigma^2}}e^{pt}dp}{\sqrt{p^2+\sigma^2}}=
J_0\left(\sigma\sqrt{t^2-x^2}\right)\Theta(t-|x|),
\end{equation}
are held on, where $\Theta(x)$ is the Heaviside function, and $J_0(x)$ is the zero-order Bessel function of real argument.

The properties-in-details of the function 
$$G_{\alpha, \gamma}(x,t)=\frac{1}{2}\int\limits_{|x|}^{\infty}
e^{i\gamma	\tau}J_0\left(\sigma\sqrt{\tau^2-x^2}\right) 
\frac{1}{t}\phi\left(-\alpha, 0;-\frac{\tau}{t^{\alpha}}\right)d\tau$$
have been analyzed in the \cite{mamchuev-2017-fcaa}.

Note, in limit of the FOP $\alpha=1,$ the function $G_{\alpha, \gamma}(x,t)$ can be cast into the form
$$\lim\limits_{\alpha\to 1}G_{\alpha, \gamma}(x,t)=
e^{i\gamma	|x|}J_0\left(\sigma\sqrt{t^2-x^2}\right)\Theta(t-|x|).$$

The estimate
\begin{equation}\label{eq210}
\left|\frac{\partial^m}{\partial x^m}D_{0t}^{\nu}G_{\alpha, \gamma}(x,t)\right|\leq C |x|^{-\theta}y^{\alpha(1-m+\theta)-\nu-1}, \quad \theta\geq 0
\end{equation}  
holds on 
for arbitrary $m\in {\mathbb N}\cup \{0\}$
and  $\nu\in {\mathbb R}$, 
where $C$ is the positive constant \cite{mamchuev-2017-fcaa}.

Accordingly, the following estimate
\begin{equation}\label{eq211}
\left|\frac{\partial^m}{\partial x^m}D_{0t}^{\nu}G_{\alpha, \gamma}(x,t)\right|\leq C \exp\left(-\sigma_0 |x|^{\varepsilon}t^{-\alpha\varepsilon}\right), 
\end{equation}  
holds on
as $|x|\to \infty$
for all $t<\infty$ \cite{mamchuev-2017-fcaa},
and for arbitrary
$m\in{\mathbb N}\cup \{0\},$ $ \nu\in {\mathbb R}$;
here and further $\sigma_0<(1-\alpha)\alpha^{\alpha\varepsilon},$ $\varepsilon=\frac{1}{1-\alpha};$  $C$ is the positive constant.

{\bf Lemma 2.1.} 
\cite[p. 172]{mam-mon-2013}. 
For $0<\mu,\nu <1$ the  estimate
\begin{equation}\label{eq212}
\int\limits_{x-\delta_1}^{x+\delta}
|x-y|^{\mu-1}|y-\xi|^{\nu-1}dy\leq C_{\gamma}|x-\xi|^{\gamma-1},
\end{equation}
holds on, where 
$$\gamma=\min\{\mu,\nu\}, \quad
C_{\gamma}=\left[\frac{2}{\gamma}+B(\mu,\nu)\right](\delta_1+\delta_2)^{\mu+\nu-\gamma}.$$

{\bf Lemma 2.2.}  
\cite[p. 176]{mam-mon-2013}. 
Let $0<\eta<s<t<T,$ $\sigma_2<\sigma_3<\sigma_1,$
then the inequality 
\begin{equation}\label{eq213}
\int\limits_{x\pm\delta}^{\pm\infty}
\exp\left[-\frac{\sigma_1 |x-y|^{\varepsilon}}{(t-s)^{\alpha\varepsilon}}\right]
\exp\left[-\frac{\sigma_2 |y-\xi|^{\varepsilon}}{(s-\eta)^{\alpha\varepsilon}}\right]
dy\leq C\exp\left[-\frac{\sigma_2 |x-\xi|^{\varepsilon}}{(t-\eta)^{\alpha\varepsilon}}\right],
\end{equation}
takes place, where  $$C=\frac{\sigma_1-\sigma_3}{\varepsilon T^{\varepsilon\alpha}}\delta^{1-\varepsilon}
\exp\left[-\frac{\sigma_1-\sigma_3}{T^{\varepsilon\alpha}}\delta^{\varepsilon}\right].$$

Furthermore, one needs the following formulae (see \cite[p. 201]{Prudnikov-1983}) 
\begin{equation}\label{eq214}
\int\limits_{0}^{\tau}
\frac{\eta}{\sqrt{\tau^2-\eta^2}}\cos\left(\rho\sqrt{\tau^2-\eta^2} \right) J_0\left(\sigma\eta\right)d\eta=\frac{1}{k}
\sin\left(k\tau\right), 
\end{equation}
\begin{equation}\label{eq215}
\int\limits_{0}^{\tau}
\frac{J_1\left(\sigma\eta\right) }{\sqrt{\tau^2-\eta^2}}\cos\left(\rho\sqrt{\tau^2-\eta^2}\right) d\eta=
\frac{1}{\sigma\tau}\cos(\rho\tau)-
\frac{1}{\sigma\tau}\cos\left(k\tau\right), 
\end{equation}
where $k=\sqrt{\sigma^2+\rho^2}.$

Under formula (\ref{eq215}), one can obtain
$$\int\limits_{0}^{\tau}
\sin\left(\rho\sqrt{\tau^2-\eta^2} \right) J_1\left(\sigma\eta\right)d\eta=
\rho\int\limits_{0}^{\tau}\left(\xi
\int\limits_{0}^{\xi}
\frac{\cos\left(\rho\sqrt{\xi^2-\eta^2} \right)}{\sqrt{\xi^2-\eta^2}} J_1\left(\sigma\eta\right)d\eta \right) d\xi=$$
\begin{equation}\label{eq216}
=\frac{\rho}{\sigma}
\int\limits_{0}^{\tau}\left[ \cos(\rho\xi)-
\cos\left(k\xi\right)\right] d\xi=
\frac{1}{\sigma}\sin(\rho\tau)-
\frac{\rho}{\sigma k}\sin(k\tau). 
\end{equation} 

The Stankovic's transformation integral (see \cite[p. 84]{pskhu-mon-2005})
\begin{equation}\label{eq217}
\int\limits_{0}^{\infty}\exp(\lambda\tau)t^{\nu-1}\phi\left(-\mu,\nu; -\tau t^{-\mu}\right)d\tau= t^{\mu+\nu-1}E_{\mu}\left(-\lambda t^{\mu}; \mu+\nu\right),
\end{equation}
takes place  for any
$\lambda\in {\mathbb C},$
$\mu\in(0,1),$
$\nu\in {\mathbb R},$
where 
$$E_{\rho}(z;\mu)=\sum\limits_{k=0}^{\infty}\frac{z^k}{\Gamma(\mu+\rho k)}$$
is the Mittag--Leffler-type function
\cite[c. 117]{Djrb-1966}. 

\section{Matrix integral formalism of the DO Cauchy problem}\label{sec:3} 
\setcounter{equation}{0}

Let us convert the DO Cauchy problem in the 'differential form' (\ref{eq11})--(\ref{eq12}) to the matrix Fredholm--Volterra-type integral equation of the second kind. The latter takes a special sense to build up the Resolvent-function solution of the DO Cauchy problem in terms of the
Liouville--Neumann-type series, which in turn is very important to computer-aided modeling
and recovering the crystal-lattice displacement field function $f({\bf R})$ from the X-ray DO
microtomography data.

The system of differential TT-type equations (\ref{eq11}) may be rewritten into the form
\begin{equation}\label{eq31}
\left (
\begin{array}{cc}
O_{-}^{\alpha}-i\gamma		& 0  	 \\
 0   	&     O_{+}^{\alpha}-i\gamma	
\end{array}
\right ) {\bf E}=i\sigma {\bf KE},
\end{equation}
where
$$O_{+}^{\alpha}=\partial_{0t}^{\alpha}+\frac{\partial}{\partial x},
\quad O_{-}^{\alpha}=\partial_{0t}^{\alpha}-\frac{\partial}{\partial x},$$
$$
{\bf E}\equiv
{\bf E}(x,t)=
\left (
\begin{array}{c}
E_0(x,t)  	 \\
E_h(x,t) 
\end{array}
\right ),
\quad
{\bf K}\equiv
{\bf K}(x,t)=
\left (
\begin{array}{cc}
0			  	&  e^{if(x,t)}  	 \\
e^{-if(x,t)}    & 0  
\end{array}
\right ). 
$$

Acting onto both sides of (\ref{eq31}) by the operator 
$diag(O_{+}^{\alpha}-i\gamma	, O_{-}^{\alpha}-i\gamma	),$ one obtains 
$$\left (
\begin{array}{cc}
O_{+}^{\alpha}-i\gamma		& 0  	 \\
 0   	&      O_{-}^{\alpha}-i\gamma	
\end{array}
\right )  
\left (
\begin{array}{cc}
O_{-}^{\alpha}-i\gamma		& 0  	 \\
 0   	&     O_{+}^{\alpha}-i\gamma	
\end{array}
\right )  
{\bf E}=$$
\begin{equation}\label{eq33}
=
i\sigma 
\left (
\begin{array}{cc}
O_{+}^{\alpha}-i\gamma		& 0  	 \\
 0   	&      O_{-}^{\alpha}-i\gamma	
\end{array}
\right )
{\bf KE}
\end{equation}
and after some straightforward routine calculations, one finds out
$$\left (\!\!
\begin{array}{cc}
O_{+}^{\alpha}-i\gamma		& 0  	 \\
 0   	&      O_{-}^{\alpha}-i\gamma	
\end{array}
\!\! \right )
{\bf KE}=
\left [\partial_{0t}^{\alpha}+
\left (\!\!
\begin{array}{cc}
1	&    0   	 \\
0   &    -1  
\end{array}
\!\!\right ) 
\frac{\partial}{\partial x}
-i\gamma \right ]({\bf KE})=
\partial_{0t}^{\alpha}({\bf KE})+$$
$$+
\left (
\begin{array}{cc}
1	&    0   	 \\
0   &    -1  
\end{array}
\right ) 
\left[
\left(\frac{\partial}{\partial x}{\bf K}\right){\bf E}+
{\bf K}\left(\frac{\partial}{\partial x}{\bf E}\right)
\right]
+{\bf K}\left(\partial_{0t}^{\alpha}{\bf E}\right)-
{\bf K}\left(\partial_{0t}^{\alpha}{\bf E}\right)
-i\gamma\,{\bf KE}=$$
$$=
{\bf K}
\left (\!\!
\begin{array}{cc}
O_{-}^{\alpha}-i\gamma	&     0 	 \\
0  	&     O_{+}^{\alpha}-i\gamma
\end{array}
\!\!\right )  
{\bf E}+\partial_{0t}^{\alpha}({\bf KE})+
\left (\!\!
\begin{array}{cc}
1	&    0   	 \\
0   &    -1  
\end{array}
\!\!\right ) 
\left(\frac{\partial}{\partial x}{\bf K}\right){\bf E}-
{\bf K}\left(\partial_{0t}^{\alpha}{\bf E}\right)=
$$
\begin{equation}\label{eq34}
={\bf K}
\left (
\begin{array}{cc}
O_{-}^{\alpha}-i\gamma	&     0 	 \\
0  	&     O_{+}^{\alpha}-i\gamma
\end{array}
\right ){\bf E}
+\left(\partial_{0t}^{\alpha}+if'_x\right)({\bf KE})
-{\bf K}\left(\partial_{0t}^{\alpha}{\bf E}\right).  
\end{equation}

The column vector ${\bf E=E}(x,t)$ is nothing else the solution of Eq. (\ref{eq31}) and ${\bf K}^2$ is equal to Unit matrix, from Eq. (\ref{eq34}) it directly follows
$$i\sigma\left (
\begin{array}{cc}
O_{+}^{\alpha}-i\gamma		& 0  	 \\
 0   	&      O_{-}^{\alpha}-i\gamma	
\end{array}
\right )
{\bf KE}=-\sigma^2 {\bf E}
+i\sigma\left(\partial_{0t}^{\alpha}+if'_x\right)({\bf KE})
-i\sigma{\bf K}\left(\partial_{0t}^{\alpha}{\bf E}\right).$$
The last allows us to rewrite down Eq. (\ref{eq33}) into the form
$$\left (
\begin{array}{cc}
(O_{+}^{\alpha}-i\gamma	)(O_{-}^{\alpha}-i\gamma	)+\sigma^2      & 0  	 \\
 0   	&      (O_{-}^{\alpha}-i\gamma	)(O_{+}^{\alpha}-i\gamma	)+\sigma^2
\end{array}
\right )  
{\bf E} =$$
\begin{equation}\label{eq35}
=i\sigma\left(\partial_{0t}^{\alpha}+if'_x\right)({\bf KE})
-i\sigma{\bf K}\left(\partial_{0t}^{\alpha}{\bf E}\right).  
\end{equation}
Keeping in mind Eq.(\ref{eq23}), and applying the double Fourier--Laplace transform to Eq. (\ref{eq35}),
one can obtain 
$$\left[(O_{+}^{\alpha}-i\gamma	)(O_{-}^{\alpha}-i\gamma	)E_0(x,t)\right]_{k,p}=
(p^{\alpha}-i\gamma	+ ik)[(O_{-}^{\alpha}-i\gamma	)E_0(x,t)]_{k,p}-$$
$$
-p^{\alpha-1}\left\{[(O_{-}^{\alpha}-i\gamma	)E_0(x,t)]_{t=0}\right\}_k=
(p^{\alpha}-i\gamma	+ ik)(p^{\alpha}-i\gamma	- ik)[E_0(x,t)]_{k,p}-$$
\begin{equation} \label{eq36}
-p^{\alpha-1}\Big\{(p^{\alpha}-i\gamma	+ ik)[E_0(x,0)]_{k}+
i\sigma [ e^{if(x,0)}E_h(x,0)]_{k}\Big\},
\end{equation}
and
$$\left[(O_{-}^{\alpha}-i\gamma	)(O_{+}^{\alpha}-i\gamma	)E_h(x,t)\right]_{k,p}=
(p^{\alpha}-i\gamma	- ik)(p^{\alpha}-i\gamma	+ ik)[E_h(x,t)]_{k,p}-$$
\begin{equation} \label{eq37}
-p^{\alpha-1}\Big\{(p^{\alpha}-i\gamma	- ik)[E_h(x,0)]_{k}+
i\sigma [ e^{-if(x,0)}E_0(x,0)]_{k}\Big\}.
\end{equation}
Here, one has used the relationships ({\it cf.} (\ref{eq31}))
$$\left[(O_{-}^{\alpha}-i\gamma	)E_0(x,t)\right]_{t=0} =i\sigma e^{if(x,0)}E_h(x,0),$$
$$\left[(O_{+}^{\alpha}-i\gamma	)E_h(x,t)\right]_{t=0}=i\sigma e^{-if(x,0)}E_0(x,0).$$
Thus from Eqs.(\ref{eq35}) -- (\ref{eq37}) it directly follows
$$\left (
\begin{array}{c}
E_0(x,t)	 \\
E_h(x,t)
\end{array}
\right )_{k,p}=
\frac{p^{\alpha-1}}{(p^{\alpha}-i\gamma	)^2+k^2+\sigma^2}\times$$
$$\times \left\{
\left (\!\!
\begin{array}{cc}
p^{\alpha}-i\gamma	+ik				& 0  	 \\
0    	& p^{\alpha}-i\gamma	-ik  
\end{array}
\!\!\right )\!\!
\left (\!\!
\begin{array}{c}
E_0(x,0)	 \\
E_h(x,0)
\end{array}
\!\!\right )_k+
i\sigma \left (\!\!
\begin{array}{c}
e^{if(x,0)}E_h(x,0)	 \\
e^{-if(x,0)}E_0(x,0)
\end{array}
\!\!\right )_k
\right\}+$$
$$+\frac{1}{(p^{\alpha}-i\gamma	)^2+k^2+\sigma^2}
\left\{
\left (
\begin{array}{cc}
\partial_{0t}^{\alpha}+if'_x				& 0  	 \\
0    	& \partial_{0t}^{\alpha}+if'_x  
\end{array}
\right )
\left (
\begin{array}{c}
i\sigma e^{if}E_h(x,t)	 \\
i\sigma e^{-if}E_0(x,t)
\end{array}
\right )-
\right. 
$$
\begin{equation}\label{eq38}
- 
\left. 
\left (
\begin{array}{cc}
0				& i\sigma e^{if}  	 \\
i\sigma e^{-if}     	& 0  
\end{array}
\right )
\left (
\begin{array}{c}
\partial_{0t}^{\alpha}E_0(x,t)	 \\
\partial_{0t}^{\alpha}E_h(x,t)
\end{array}
\right )
\right\}_{k,p}.  
\end{equation}
Applying the Efros's theorem  for operational calculus \cite[p. 512]{Lavr},
the equalities  (\ref{eq24})--(\ref{eq29}),  
and using the inverse double Fourier--Laplace transform, the following elations take
place
$$\left(\frac{1}{(p^{\alpha}-i\gamma	)^2+k^2+\sigma^2}\right)_{x,t}=\frac{1}{i(2\pi)^2}\int\limits_{-i\infty}^{i\infty}dp\int\limits_{-\infty}^{\infty}
\frac{e^{pt+ikx}}{(p^{\alpha}-i\gamma	)^2+k^2+\sigma^2}dk=$$
$$=\frac{1}{2}\int\limits_{|x|}^{\infty}
e^{i\gamma	\tau}J_0\left(\sigma\sqrt{\tau^2-x^2}\right) 
\frac{1}{t}\phi\left(-\alpha, 0;-\frac{\tau}{t^{\alpha}}\right)d\tau
=G_{\alpha, \gamma	}(x,t),$$

$$\left(\frac{p^{\alpha-1}}{(p^{\alpha}-i\gamma	)^2+k^2+\sigma^2}\right)_{x,t}=\frac{1}{i(2\pi)^2}\int\limits_{-i\infty}^{i\infty}dp
\int\limits_{-\infty}^{\infty}
\frac{p^{\alpha-1}e^{pt+ikx}}{(p^{\alpha}-i\gamma	)^2+k^2+\sigma^2}dk=
$$
$$=D_{0t}^{\alpha-1}G_{\alpha, \gamma}(x,t)=
\frac{1}{2}\int\limits_{|x|}^{\infty}
e^{i\gamma	\tau}J_0\left(\sigma\sqrt{\tau^2-x^2}\right) 
t^{-\alpha}\phi\left(-\alpha, 1-\alpha;-\frac{\tau}{t^{\alpha}}\right)d\tau.$$

As a result, from (\ref{eq38}) one obtains the DO integral matrix equation
\begin{equation}\label{eq39}
{\bf E}(x,t)=({\bf A}^{\alpha, \gamma}{\bf E}(x,t))(x,t)+({\bf B}^{\alpha, \gamma}{\bf E}(x,0))(x,t),
\end{equation}
where the following notations are introduced
$$({\bf A}^{\alpha, \gamma}{\bf E}(x,t))(x,t)=-i\sigma\int\limits_{0}^{t}dv\int\limits_{-\infty}^{\infty}
G_{\alpha, \gamma}(x-u,t-v)\cdot
\left\{D_{0v}^{\alpha-1}\frac{\partial}{\partial v}
[{\bf K}(u,v){\bf E}(u,v)]+\right.$$
$$
+\left.
\left (
\begin{array}{cc}
1		& 0  	 \\
0    	& -1  
\end{array}
\right )\left[\frac{\partial}{\partial u}{\bf K}(u,v)\right]
{\bf E}(u,v)-
{\bf K}(u,v)D_{0v}^{\alpha-1}\frac{\partial}{\partial v}
{\bf E}(u,v)\right\}du,
$$

$$({\bf B}^{\alpha, \gamma}{\bf E}(x,0))(x,t)=\int\limits_{0}^{t}dv\int\limits_{-\infty}^{\infty}
D_{tv}^{\alpha-1}G_{\alpha, \gamma}(x-u,t-v)\times$$
$$
\times \left (
\begin{array}{cc}
O^{\alpha}_+	-i\gamma				& i\sigma e^{if(u,0)}  	 \\
i\sigma e^{-if(u,0)}    	& O^{\alpha}_- -i\gamma	 
\end{array}
\right )
{\bf E}(u,0)\delta(v)du,
$$
where $\delta(v)$ is the Dirac delta-function.

Using the formula for integration by parts and its fractional analogue (\ref{eq22}), the estimates (\ref{eq210}) and (\ref{eq211}), definition of the fractional Riemann--Liouville derivative, 
one finds out
$$({\bf A}^{\alpha, \gamma}{\bf E}(x,t))(x,t)=-i\sigma\int\limits_{0}^{t}dv\int\limits_{-\infty}^{\infty}
{\bf K}_1(x,t;u,v){\bf E}(u,v)du-
$$
\begin{equation}\label{eq312}
-i\sigma
\int\limits_{-\infty}^{\infty}D_{0t}^{\alpha-1}
G_{\alpha, \gamma}(x-u,t){\bf K}(u,0){\bf E}(u,0)du,
\end{equation}
where
$$
{\bf K}_1(x,t;u,v)=
D_{vt}^{\alpha}G_{\alpha, \gamma}(x-u,t-v)\cdot{\bf K}(u,v)+$$$$+if'_{u}(u,v)\cdot G_{\alpha, \gamma}(x-u,t-v){\bf K}(u,v)
+D_{vt}^{\alpha}\left[G_{\alpha, \gamma}(x-u,t-v){\bf K}(u,v)\right].$$

In view of (\ref{eq21}) the $({\bf B}^{\alpha, \gamma}{\bf E}(x,0))(x,t)$ may be rewritten as
$$({\bf B}^{\alpha, \gamma}{\bf E}(x,0))(x,t)=
-\int\limits_{-\infty}^{\infty}D_{0t}^{2\alpha-1}
G_{\alpha, \gamma}(x-u,t){\bf E}(u,0)du-$$
$$-i\gamma	\int\limits_{-\infty}^{\infty}D_{0t}^{\alpha-1}
G_{\alpha, \gamma}(x-u,t){\bf E}(u,0)du+$$
$$+\left (
\begin{array}{cc}
1		& 0  	 \\
0    	& -1  
\end{array}
\right )
\int\limits_{-\infty}^{\infty}D_{0t}^{\alpha-1}
G_{\alpha, \gamma}(x-u,t){\bf E}'(u,0)du+$$
\begin{equation}\label{eq313}
+i\sigma
\int\limits_{-\infty}^{\infty}D_{0t}^{\alpha-1}
G_{\alpha, \gamma}(x-u,t){\bf K}(u,0){\bf E}(u,0)du.
\end{equation}

Taking into account Eqs. (\ref{eq312}), (\ref{eq313}), the integral matrix equation (\ref{eq39}) may  be to reduce  
\begin{equation}\label{eq314}
{\bf E}(x,t)+i\sigma\int\limits_{0}^{t}dv\int\limits_{-\infty}^{\infty}
{\bf K}_1(x,t;u,v){\bf E}(u,v)du=
{\bf F}(x,t),
\end{equation}
where the following notations are introduced
$$
{\bf F}(x,t)=({\bf B}^{\alpha, \gamma}{\bf E}(x,0))(x,t)-i\sigma
\int\limits_{-\infty}^{\infty}D_{0t}^{\alpha-1}
G_{\alpha, \gamma}(x-u,t){\bf K}(u,0){\bf E}(u,0)du.$$

\section{Matrix Resolvent-function solution of the DO Cauchy problem}
\label{sec:4}
\setcounter{equation}{0}


Let $A(x)$ be the matrix with entries $a_{ij}(x).$ 
By notation $|A(x)|_*$ one denotes a scalar
function taking the maximum absolute value of entries $a_{ij}(x)$ of the matrix $A(x)$ for each $x$; i.e.,
$|A(x)|_* = \max\limits_{i,j}|a_{ij}(x)|.$

Let $f(x,t)$ and $f_t(x,t)$ be the continuous, bounded functions.

From the Eqs. (\ref{eq210}) -- (\ref{eq213}) one obtains following estimates for the matrix kernel function ${\bf K}_1$
\begin{equation}\label{eq41}
|{\bf K}_1(x,t;\xi,\eta)|_*\leq C |x-\xi|^{-\theta}(t-\eta)^{\beta-1}, 
\end{equation}
\begin{equation}\label{eq42}
|{\bf K}_1(x,t;\xi,\eta)|_*\leq C(t-\eta)^{\beta-1}
\exp\left[-\frac{\sigma_0 |x-\xi|^{\varepsilon}}{(t-\eta)^{\alpha\varepsilon}}\right], 
\end{equation}
where $\beta=\alpha\theta,$ $\theta\in(0,1),$
$\sigma_0<(1-\alpha)\alpha^{\alpha\varepsilon},$
$\varepsilon=\frac{1}{1-\alpha}.$

Further, let us find out the corresponding estimates for iterative kernels functions, namely:
$${\bf K}_n(x,t;\xi,\eta)=\int\limits_{\eta}^{t}dv\int\limits_{-\infty}^{\infty}
{\bf K}_{n-1}(x,t;u,v){\bf K}_1(u,v;\xi,\eta)du.$$

From estimates (\ref{eq41}) and (\ref{eq42})  it follows
$$|{\bf K}_2(x,t;\xi,\eta)|_*\leq $$
$$ \leq  
\int\limits_{\eta}^{t}dv\left(\int\limits_{x-\delta_1}^{x+\delta_2}+\int\limits_{-\infty}^{x-\delta_1}+\int\limits_{x+\delta_2}^{\infty}\right)
|{\bf K}_1(x,t;u,v){\bf K}_1(u,v;\xi,\eta)|_*du\leq
$$
\begin{equation}\label{eq43}
\leq C C_{2,\beta,\theta}
(t-\eta)^{2\beta -1}\left(|x-\xi|^{-\theta}+
\exp\left[-\frac{\sigma_0 |x-\xi|^{\varepsilon}}{(t-\eta)^{\alpha\varepsilon}}\right]\right),
\end{equation}
where 
$$C_{2,\beta,\theta}=\frac{\Gamma^2(\beta)C_{1-\theta}}{\Gamma(2\beta)}.$$
Accordingly, from (\ref{eq43}), keeping in mind $0<e^{-s}<1,$ for $s>0,$ one obtains
\begin{equation}\label{eq44}
|{\bf K}_2(x,t;\xi,\eta)|_*\leq C 
C_{2,\beta,\theta}
|x-\xi|^{-\theta}(t-\eta)^{2\beta -1}, 
\end{equation}
\begin{equation}\label{eq45}
|{\bf K}_2(x,t;\xi,\eta)|_*\leq C C_{2,\beta,\theta}
(t-\eta)^{2\beta -1}
\exp\left[-\frac{\sigma_0 |x-\xi|^{\varepsilon}}{(t-\eta)^{\alpha\varepsilon}}\right]. 
\end{equation}
Under the above, the following estimates take place
\begin{equation}\label{eq46}
|{\bf K}_n(x,t;\xi,\eta)|_*\leq C 
C_{n,\beta,\theta} {\Gamma(n\beta)}|x-\xi|^{-\theta}(t-\eta)^{n\beta -1}, 
\end{equation}
\begin{equation}\label{eq47}
|{\bf K}_n(x,t;\xi,\eta)|_*\leq C 
C_{n,\beta,\theta}(t-\eta)^{n\beta -1}
\exp\left[-\frac{\sigma_0 |x-\xi|^{\varepsilon}}{(t-\eta)^{\alpha\varepsilon}}\right].
\end{equation}
where 
$$C_{n,\beta,\theta}=\frac{\Gamma^n(\beta)C_{1-\theta}^{n-1}}{\Gamma(n\beta)}.$$
The matrix Resolvent-function solution of the Eq. (\ref{eq314}) takes the form
\begin{equation}\label{eq48}
{\bf R}(x,t;\xi,\eta)=\sum\limits_{n=1}^{\infty}{\bf K}_{n}(x,t;\xi,\eta).
\end{equation}

From estimates (\ref{eq41})--(\ref{eq47}), it follows the Liouville--Neumann-type series (\ref{eq48}) converges and can be estimated as a value of
$$|{\bf R}(x,t;\xi,\eta)|_*\leq C \Gamma(\beta) (t-\eta)^{\beta-1}
E_{\beta}\left[C_{1-\theta}\Gamma(\beta)(t-\eta)^{\beta};\beta\right]\times$$
$$
\times
\left(|x-\xi|^{-\theta}+\exp\left[-\frac{\sigma_0 |x-\xi|^{\varepsilon}}{(t-\eta)^{\alpha\varepsilon}}\right]\right).
$$

On the other hand, since the Mittag-Leffler-type function $E_{\beta}(s;\beta)$ is bounded on the interval $[0,T]$, the estimates of (\ref{eq41})--(\ref{eq42}) are proved and valid also for the matrix Resolvent-function solution as a whole.

So, one can state the unique solution of the matrix integral Eq. (\ref{eq314}) takes the form
$${\bf E}(x,t)= {\bf F}(x,t)+\int\limits_{0}^{t}\int\limits_{-\infty}^{+\infty}{\bf R}(x,t;\xi,\eta) {\bf F}(\xi,\eta)d\xi d\eta.$$

\section{The DO Cauchy problem. The FOP $\alpha=1$}
\label{sec:5}
\setcounter{equation}{0}

In the case when FOP  $\alpha=1,$ the basic system (\ref{eq31}) reduces to
$$\left (
\begin{array}{cc}
\frac{\partial}{\partial t}-\frac{\partial}{\partial x}	& 0  	 \\
 0   	&      \frac{\partial}{\partial t}+\frac{\partial}{\partial x}
\end{array}
\right )  
\left (
\begin{array}{c}
E_0(x,t)	 \\
E_h(x,t)
\end{array}
\right ) =$$
$$
=i
\left (
\begin{array}{cc}
\gamma						& \sigma \exp[if(x,t)]  	 \\
\sigma\exp[-if(x,t)]    	& \gamma	     
\end{array}
\right )  
\left (
\begin{array}{c}
E_0(x,t)	 \\
E_h(x,t)
\end{array}
\right ).
$$

As was shown earlier \cite{MamChukh-2023}, one can find out that for $\alpha\to 1$ integral equation (\ref{eq39}) can be written down into the form
$$
{\bf E}(x,t)=({\bf A}^{1,\gamma}{\bf E}(x,t))(x,t)+({\bf B}^{1,\gamma}{\bf E}(x,0))(x,t),
$$
where
$$({\bf A}^{1,\gamma}{\bf E}(x,t))(x,t)=i\sigma
\int\limits_0^t dv
\int\limits_{x-t+v}^{x+t-v}
G_{1,\gamma}(x-u,t-v)\times $$
$$\times
\left(
\begin{array}{cc}
0									&	{O_+\, e^{if(u,v)}} \\
{O_-\, e^{-if(u,v)}}	&		0
\end{array}
\right){\bf E}(u,v)du,$$
$$({\bf B}^{1,\gamma}{\bf E}(x,0))(x,t)=
\int\limits_0^t dv
\int\limits_{x-t+v}^{x+t-v}
G_{1,\gamma}(x-u,t-v)\times $$
$$\times
\left(
\begin{array}{cc}
O_+	-i\gamma				&	i\sigma e^{if(u,0)} \\
i\sigma e^{-if(u,0)}	&		O_- -i\gamma 
\end{array}
\right)\delta(v){\bf E}(u,0)du,$$
and respectively,
$$
G_{1,\gamma}(x,t)=\frac{1}{2}e^{i\gamma t}J_0\left(\sigma \sqrt{t^2-x^2}\right)\Theta[t-|x|].
$$

\section{The DO Cauchy problem. The crystal-lattice displacement field function $f({\bf R})=a x +b$}
\label{sec:6}
\setcounter{equation}{0}

Here we will build up a solution of the basic fractional DO Cauchy problem when the crystal-lattice displacement field function 
$f({\bf R})=a x+b$.

After trivial exponentional substitutions for the wave amplitudes $E_0(x, t),$ $E_h(x, t),$ the system (\ref{eq11}) can written down as (for simplicity, further, the same notations for the wave amplitudes $E_0(x, t),$ $E_h(x, t),$ are to be saved)
\begin{equation}\label{eq61}
\begin{split}
\left (\partial_{0t}^{\alpha}-\frac{\partial}{\partial x}\right)E_0(x,t)=i\gamma	 E_0(x,t)+i\sigma e^{ia x}E_h(x,t), \,\,\,\, \\
\left (\partial_{0t}^{\alpha}+\frac{\partial}{\partial x}
\right )E_h(x,t)=i\sigma e^{-ia x}E_0(x,t)+i\gamma	 E_h(x,t).
\end{split}
\end{equation}
Substituting the functions $E_0(x, t),$ $E_h(x, t)$ as
$$E_0(x,t)=\exp\left(i\frac{a x}{2}\right){\mathcal E}_0(x,t),
\quad
E_h(x,t)=\exp\left(-i\frac{a x}{2}\right){\mathcal E}_h(x,t),$$
one obtains 
\begin{equation}\label{eq62}
\begin{split}
  \left (\partial_{0t}^{\alpha}-\frac{\partial}{\partial x}\right){\mathcal E}_0(x,t)=i\left(\gamma+\frac{a}{2}\right) {\mathcal E}_0(x,t)+i\sigma {\mathcal E}_h(x,t), \\
  \left (\partial_{0t}^{\alpha}+\frac{\partial}{\partial x}
\right ){\mathcal E}_h(x,t)=i\sigma {\mathcal E}_0(x,t)+
i\left(\gamma+\frac{a}{2}\right) {\mathcal E}_h(x,t),
\end{split}
\end{equation}
and the initial condition
\begin{equation}\label{eq63}
\left (
\begin{array}{c}
{\mathcal E}_0(x,0)	 \\
{\mathcal E}_h(x,0)
\end{array}
\right ) =
\left (
\begin{array}{c}
\psi_0(x)	 \\
\psi_h(x)
\end{array}
\right )=
\left (
\begin{array}{c}
e^{-ia x/2}\varphi_0(x)	 \\
e^{ia x/2}\varphi_h(x)
\end{array}
\right ), 
\quad -\infty<x<\infty.
\end{equation}

Then, let us introduce the notations
\begin{equation}\label{eq64}
{\bf \Gamma}(x,t)=\frac{1}{2}\int\limits_{|x|}^{\infty}
g(t,\tau){\bf Q}(x,\tau)d\tau+{\bf \Gamma}_0(x,t),
\end{equation}
\begin{equation}\label{eq65}
{\bf Q}(x,\tau)=
\left [
\begin{array}{cc}
-\sigma\frac{\tau-x}{\sqrt{\tau^2-x^2}}J_1(\sigma\sqrt{\tau^2-x^2})	& i\sigma J_0(\sigma\sqrt{\tau^2-x^2}) \\
i\sigma J_0(\sigma\sqrt{\tau^2-x^2}) &			-\sigma\frac{\tau+x}{\sqrt{\tau^2-x^2}}J_1(\sigma\sqrt{\tau^2-x^2})
\end{array}
\right ], 
\end{equation}

\begin{equation}\label{eq66}
{\bf \Gamma}_0(x,t)=g(t,|x|)
\left [
\begin{array}{cc}
\Theta(-x)	& 0 \\
0		&  \Theta(x)
\end{array}
\right ], 
\end{equation}
$$g(t,\tau)=
\frac{e^{i(\gamma+a/2)\tau}}{t}\phi\left(-\alpha,0; -\tau t^{-\alpha}\right),$$
$\Theta(x)$ is the Heaviside function,
$J_m(z)$ is the $m$-order Bessel function of the argument $z$. 

According to \cite{mamchuev-2021-aman}, the basic matrix solution of the DO Cauchy problem (\ref{eq62}) -- (\ref{eq63}) has the form
\begin{equation}\label{eq67}
{\mathcal E}(x,t)=\int\limits_{-\infty}^{\infty}D_{0t}^{\alpha-1}\Gamma(x-\xi,t)\psi(\xi)d\xi,
\end{equation}
in the class of function  
$${\mathcal E}(x,t)\in C(\overline{\Omega}), \quad \partial_{0t}^{\alpha}{\mathcal E}(x,t), \frac{\partial}{\partial x} {\mathcal E}(x,t)\in C(\Omega),$$ 
where  $\psi(x)=[\psi_0(x), \psi_h(x)]^{tr}\in C(-\infty,\infty)$ is the function, which
satisfies the Hölder condition, and the following relation as $|x|\to\infty$ 
$$\psi(x)=O(\exp(\rho|x|^{\varepsilon})) \quad 
\varepsilon=\frac{1}{1-\alpha}, \quad \rho<(1-\alpha)(\alpha T^{-1})^{\frac{\alpha}{1-\alpha}}.$$ 

Underline the solution to the Cauchy problem (\ref{eq62}) -- (\ref{eq63}) is unique in the class of functions, which satisfy the condition 
$${\mathcal E}(x,t)=O(\exp(k|x|^{\varepsilon})), \quad \text{ as} 
\quad |x|\to\infty,$$
for some $k>0.$

From Eqs. (\ref{eq63}) -- (\ref{eq67}) it follows  the solution of the Cauchy problem (\ref{eq61}), (\ref{eq12})
can be cast into the form 
\begin{equation}\label{eq68}
{\bf E}(x,t)=\int\limits_{-\infty}^{\infty}D_{0t}^{\alpha-1}{\bf G}(x,\xi,t)\varphi(\xi)d\xi,
\end{equation}
where the following notations are introduced
\begin{equation}\label{eq69}
{\bf G}(x,\xi,t)=\frac{1}{2}\int\limits_{|x-\xi|}^{\infty}
g(t,\tau){\bf S}(x,\xi,\tau)d\tau+{\bf G}_0(x,\xi,t),
\end{equation}
$$
{\bf S}(x,\xi,\tau)=
\left [\!\!
\begin{array}{cc}
-\sigma e^{iaX_1/2}\frac{\tau-X_1}{\sqrt{\tau^2-X_1^2}}h_1(\tau, X_1)	& 
i\sigma e^{iaX_2/2}h_0(\tau, X_1) \\
i\sigma e^{-iaX_2/2}h_0(\tau, X_1) &			
-\sigma e^{-iaX_1/2}\frac{\tau+X_1}{\sqrt{\tau^2-X_1^2}}h_1(\tau, X_1)
\end{array}
\!\!\right ], 
$$
$$h_0(\tau, X_1)=J_0(\sigma\sqrt{\tau^2-X_1^2}),
\quad 
h_1(\tau, X_1)=J_1(\sigma\sqrt{\tau^2-X_1^2})$$
\begin{equation}\label{eq610}
{\bf G}_0(x,\xi,t)=\left [\!\!
\begin{array}{cc}
e^{-i\gamma X_1}\Theta(-X_1)	& 0 \\
0							& e^{i\gamma X_1}\Theta(X_1)
\end{array}
\!\!\right ]
\frac{1}{t}\phi\left(-\alpha,0; -|X_1| t^{-\alpha}\right), 
\end{equation}
$$X_1=x-\xi, \quad X_2=x+\xi.$$

\section{Case of the function  $f({\bf R})=a x+b$ and $E_0(x, 0) = 1,$ $E_h(x, 0)=0$}
\label{sec:7}
\setcounter{equation}{0}

Let us consider the case when the initial conditions for the wavefield amplitudes $E_0(x,0)$ and $E_h(x,0)$ are constant and  equal to
\begin{equation}\label{eq71}
E_0(x,0)=\varphi_0(x)\equiv 1, \quad
E_h(x,0)=\varphi_h(x)\equiv 0.
\end{equation}
Then, the formulae (\ref{eq68})--(\ref{eq610}) can be simplified and expressed in terms of the Mittag--Leffler-type functions.

Keeping in mind Eq. (\ref{eq71}), Eq. (\ref{eq68}) for ${\bf E}(x, t)$ can be written down as
$${\bf E}(x,t)=\int\limits_{-\infty}^{\infty}D_{0t}^{\alpha-1}{\bf G}(x,\xi,t)
\left (
\begin{array}{c}
1	 \\
0
\end{array}
\right )
d\xi=$$
$$=\frac{1}{2}\int\limits_{-\infty}^{\infty}d\xi \int\limits_{|x-\xi|}^{\infty}
D_{0t}^{\alpha-1}g(t,\tau){\bf S}(x,\xi,\tau)
\left (
\begin{array}{c}
1	 \\
0
\end{array}
\right )d\tau+
$$
\begin{equation}\label{eq72}
+\int\limits_{-\infty}^{\infty}D_{0t}^{\alpha-1}{\bf G}_0(x,\xi,t)\left (
\begin{array}{c}
1	 \\
0
\end{array}
\right )d\xi
\equiv {\bf I}_1(x,t)+{\bf I}_2(x,t).
\end{equation}
By changing the integration order, let us evaluate the integral ${\bf I}_1(x,t)$ 
$${\bf I}_1(x,t)=\frac{1}{2}\int\limits_{0}^{\infty}D_{0t}^{\alpha-1}g(t,\tau)d\tau
\int\limits_{x-\tau}^{x+\tau}{\bf S}(x,\xi,\tau)
\left (
\begin{array}{c}
1	 \\
0
\end{array}
\right )d\xi=$$
\begin{equation}\label{eq73}
=\frac{1}{2}\int\limits_{0}^{\infty}D_{0t}^{\alpha-1}g(t,\tau)d\tau
\int\limits_{-\tau}^{\tau}{\bf S}(x,x-\eta,\tau)
\left (
\begin{array}{c}
1	 \\
0
\end{array}
\right )d\eta.
\end{equation}
From the Eqs. (\ref{eq214})--(\ref{eq216}), it directly follows up
\begin{equation}\label{eq74}
\int\limits_{-\tau}^{\tau}e^{\pm ia\eta/2}
J_0(\sigma\sqrt{\tau^2-\eta^2})d\eta
=\frac{2}{k} \sin (k\tau),
\end{equation} 
$$\int\limits_{-\tau}^{\tau}e^{\pm ia\eta/2}
\frac{\tau\mp\eta}{\sqrt{\tau^2-\eta^2}}J_1(\sigma\sqrt{\tau^2-\eta^2})d\eta=$$
\begin{equation}\label{eq75}
=-\frac{1}{\sigma}\left[2 \cos(k\tau)- i\frac{a}{k}\sin (k\tau)-2e^{-ia\tau/2}\right],
\end{equation}
where $k=\sqrt{a^2/4+\sigma^2}.$

Accordingly, from Eqs. (\ref{eq74}), (\ref{eq75}) one obtains
$$
\int\limits_{-\tau}^{\tau}S_{11,22}(x,x-\eta,\tau)d\eta=
-\sigma\int\limits_{-\tau}^{\tau}e^{\pm ia\eta/2}
\frac{\tau\mp\eta}{\sqrt{\tau^2-\eta^2}}J_1(\sigma\sqrt{\tau^2-\eta^2})d\eta=$$
\begin{equation}\label{eq76}
=2 \cos(k\tau)- i\frac{a}{k}\sin (k\tau)-2e^{-ia\tau/2},
\end{equation}
$$
\int\limits_{-\tau}^{\tau}S_{12,21}(x,x-\eta,\tau)d\eta=
i\sigma e^{\pm i ax}\int\limits_{-\tau}^{\tau}e^{\mp ia\eta/2}
J_0(\sigma\sqrt{\tau^2-\eta^2})d\eta=$$
\begin{equation}\label{eq77}
=i\sigma \frac{2}{k} e^{\pm i ax}\sin (k\tau),
\end{equation}
where $S_{ij}(x,x-\eta,\tau)$ $(i, j=1,2)$ are the elements of matrix ${\bf S}(x,x-\eta,\tau).$

Following equalities
$$e^{ia_1\tau}\left[2 \cos(k\tau)- i\frac{a}{k}\sin (k\tau)\right]=$$
\begin{equation}\label{eq78}
=\left(1- \frac{a}{2k}\right)e^{i(a_1+k)\tau}+
\left(1+ \frac{a}{2k}\right)e^{i(a_1-k)\tau},
\end{equation}
\begin{equation}\label{eq79}
-2e^{i a_1\tau}e^{i a\tau/2}=-2e^{i \gamma\tau},
\end{equation}
\begin{equation}\label{eq710}
e^{i a_1\tau}\sin (k\tau)=\frac{i}{2}e^{i(a_1-k)\tau}
-\frac{i}{2}e^{i(a_1+k)\tau},
\end{equation}
hold, where 
$$a_1=\gamma+\frac{a}{2},  \quad k=\sqrt{a^2/4+\sigma^2}.$$
Next, exploiting calculations (\ref{eq76})--(\ref{eq710}), from (\ref{eq73}) one obtains
\begin{equation}\label{eq711}
{\bf I}_1(x,t)=\int\limits_{0}^{\infty} t^{-\alpha}
\phi\left(-\alpha, 1-\alpha; -\tau t^{-\alpha}\right)
{\bf N}(x,\tau)\left (
\begin{array}{c}
1	 \\
0	
\end{array}
\right )
d\tau,
\end{equation}
where
$${\bf N}(x,\tau)={\bf N}_1(x)e^{i(a_1+k)\tau}	
+{\bf N}_2(x)e^{i(a_1-k)\tau}	
-
\left (
\begin{array}{cc}
1	& 0 \\
0	& 1
\end{array}
\right )
e^{i\gamma\tau},$$
$${\bf N}_{1,2}(x)=\frac{1}{2}\left (
\begin{array}{ll}
1\mp \frac{a}{2k}&  
\pm\frac{\sigma}{k} e^{iax} \\
\pm\frac{\sigma}{k} e^{-iax}	&  
1\mp\frac{a}{2k}
\end{array}
\right ).$$
From (\ref{eq610}) one finds out
\begin{equation}\label{eq712}
{\bf I}_2(x,t)=
\int\limits_{0}^{\infty} t^{-\alpha}
\phi\left(-\alpha, 1-\alpha; -\eta t^{-\alpha}\right)e^{i\gamma\eta}
\left (
\begin{array}{c}
1	 \\
0	
\end{array}
\right )d\eta.
\end{equation}

Applying formula (\ref{eq217}) to equalities (\ref{eq711}) and (\ref{eq712}), the total solution (\ref{eq72}) can be cast into the form
\begin{equation}\label{eq713}
{\bf E}(x,t)=\frac{1}{2}
\left (
\begin{array}{lr}
1-\frac{a}{2k}	&  
1+\frac{a}{2k}\\
\frac{\sigma}{k} e^{-iax}	&  
-\frac{\sigma}{k} e^{-iax}
\end{array}
\right )
\left (
\begin{array}{c}
E_{\alpha, 1}\left(i(a_1+k)t^{\alpha}\right)	 \\
E_{\alpha, 1}\left(i(a_1-k)t^{\alpha}\right)
\end{array}
\right ).
\end{equation}

In the case when the FOP $\alpha\to 1$, the total solution (\ref{eq713}) is reduced to
\begin{equation}\label{eq714}
{\bf E}(x,t)=\frac{1}{2}
\left (
\begin{array}{lr}
1-\frac{a}{2k}	&  
1+\frac{a}{2k}\\
\frac{\sigma}{k} e^{-iax}	&  
-\frac{\sigma}{k} e^{-iax}
\end{array}
\right )
\left (
\begin{array}{c}
e^{i(a_1+k)t}	 \\
e^{i(a_1-k)t}
\end{array}
\right ).
\end{equation}
Finally, putting on the parameters $\gamma=0$ and $a=0,$  the solution (\ref{eq714}) goes into the
pendulum solutions of the canonic TT equations regarding the Bloch wavefield
amplitudes in the perfect crystal (cf.  \cite{MamChukh-2023})
$$E_0=\cos \sigma t, \quad E_h=i\sin \sigma t.$$

\section{Conclusion}

In this paper, a goal of our study is to elaborate the novel mathematics framework for
processing the DO data based on a concept of the matrix integral equation to solve the inverse
Radon problem that arose in the computer microtomography (see, e.g., \cite{Authier-2001}, \cite{Chukhovskii-2020}).

In the case when FOP $\alpha=1$, the above results obtained can be applied to the
standard X-ray DO using the computer diffraction microtjmography technique (cf. \cite{MamChukh-2023}).

In contrast to \cite{Honkanen-2018}, \cite{Chukhovskii-2020}, some advantage of our study is to open a window to develop the
theory models of the X-ray and electron diffraction scattering by imperfect crystals. The
fractional DO theory operates based on the matrix Fredholm--Volterra integral equation of the second kind.
The matrix Resolvent-function solution of the
Fredholm--Volterra integral equation of the second kind allows us to model the X-ray
diffraction scattering by imperfect crystals. Besides, it can modify the optimizing
procedure of the $\chi^2$-target function in trial to solve the inverse Radon problem in the
DO microtomography investigating the crystal-lattice defects structure.

The common approach of the DO imaging, the FOP $\alpha=1$, has to be rethinked using the
fractional derivatives formalism $\partial_{0t}^{\alpha}E$ in the sense of the Gerasimov--Caputo approach (\ref{eqgc}).
Accordingly, in the fractional DO, it is also worth noting alternative interest due to the
possibility of using the fractional Riemann–Liouville derivative $D_{0t}^{\alpha}E$  (\ref{eqrl}). The latter is a
good topic for future work in the fractional DO.




 \bigskip \smallskip

 \it
 
  \noindent
$^1$ Institute of Applied Mathematics and Automation, \\
Kabardino-Balkarian Scientific Center RAS,\\
``Shortanov" Str., 89A \\
360000  Nal'chik, Russian Federation \\[4pt]
e-mail: mamchuev@rambler.ru (Corr. author)\\[4pt]
e-mail: f\_chukhov@yahoo.ca \\[12pt]
\hfill Received: March 20, 2024 \\[12pt]

\end{document}